\documentclass [11pt]{article}
\usepackage{amsthm}
\usepackage{amsmath}
\usepackage{amssymb}
\usepackage{enumerate}
\usepackage{epsfig}
\usepackage{srcltx}
\usepackage[dvipsnames,usenames]{color}

\newcounter{contador}
\newcounter{contador2}

\newtheorem{propo}[contador]{Proposition}
\newtheorem{teo}[contador]{Theorem}
\newtheorem{lem}[contador]{Lemma}
\newtheorem{nota}[contador]{Remark}

\newtheorem{corol}[contador]{Corollary}

\newcommand{\rec}{\noindent}    

\newcommand{\Qm}{{Q^+}}

\newcommand{\R}{{\mathbb R}}

\newcommand{\N}{{\mathbb N}}
\newcommand{\Q}{{\mathbb Q}}

\newcommand{\U}{{\cal{U}}}

\title{Integrability and non-integrability\\ of periodic
non-autonomous Lyness recurrences\footnote{{\bf Acknowledgements}.
GSD-UAB and CoDALab Groups are supported by the Government of
Catalonia through the SGR program. They are also supported by MCYT
through grants MTM2008-03437 (first and second authors);
DPI2008-06699-C02-02 and  DPI2011-25822 (third author).}}

\author{Anna Cima$^{(1)}$, Armengol Gasull$^{(1)}$ and V\'{\i}ctor Ma\~{n}osa $^{(2)}$
  \\*[.1truecm]
{\small \textsl{$^{(1)}$ Dept. de Matem\`{a}tiques, Facultat de
Ci\`{e}ncies,}}
\\*[-.25truecm] {\small \textsl{Universitat Aut\`{o}noma de Barcelona,}}
\\*[-.25truecm] {\small \textsl{08193 Bellaterra, Barcelona, Spain}}
\\*[-.25truecm] {\small \textsl{cima@mat.uab.cat, gasull@mat.uab.cat}}
\\*[-.25truecm]
\\*[-.25truecm] {\small \textsl{$^{(2)}$ Dept. de Matem\`{a}tica Aplicada III (MA3),}}
\\*[-.25truecm] {\small \textsl{Control, Dynamics and Applications Group (CoDALab)}}
\\*[-.25truecm] {\small \textsl{Universitat Polit\`{e}cnica de Catalunya (UPC)}}
\\*[-.25truecm] {\small \textsl{Colom 1, 08222 Terrassa, Spain}}
\\*[-.25truecm] {\small \textsl{victor.manosa@upc.edu}}}

\begin{document}
\maketitle

\begin{abstract}
This paper studies non-autonomous Lyness type recurrences of the
form $x_{n+2}=(a_n+x_{n+1})/x_{n}$, where $\{a_n\}$ is a
$k$-periodic sequence of positive numbers with primitive period $k$.
We show that for the cases $k\in\{1,2,3,6\}$ the behavior of the
sequence $\{x_n\}$ is simple (integrable) while for  the remaining
cases satisfying this behavior can be much more complicated
(chaotic). We also show that the cases where $k$ is a multiple of 5
present some different features.
\end{abstract}

\rec {\sl 2000 Mathematics Subject Classification:} \texttt{37C55,
39A11, 39A20}

\rec {\sl Keywords:}  Integrability and non-inte\-gra\-bility of
discrete systems, numerical chaos, periodic difference equations,
QRT maps, rational and meromorphic first integrals.\newline

\section{Introduction and main results}

 The dynamical study of the Lyness difference equation (\cite{BR,BC,ER,Z}) and
its generalizations to higher order Lyness type equations
(\cite{BR3,CL, CGM07-1, Gao}), or to difference equations with
periodic coefficients (\cite{CGM09,FJL,JKN,KN,RGW}), has been the
focus of an active research activity in the last two decades. In
more recent dates, Lyness type equations have also been approached
using different points of view: from algebraic geometry
(\cite{B,D,U}) to the theory of discrete integrable systems
(\cite{CGMJPhA,GRT,GRTW,RGW,TvQ}).

This paper deals with the problem of the integrability and
non-integrability of non-autonomous planar Lyness difference
equations of the form
\begin{equation}\label{eq}
x_{n+2}\,=\,\frac{a_n+x_{n+1}}{x_n},
\end{equation}
where $\{a_n\}$ is a cycle of $k$ positive numbers, i.e. $a_{n+k}=a_n$ for all
$n\in\N,$  being $k$ the primitive period and we consider positive initial
conditions $x_1$ and $x_2$. As we will see, the behavior of the sequences
$\{x_n\}$ can be  essentially different according to  whether $k\in\{1,2,3,6\}$,
$k$ is a multiple of $5$ or it is not.

 In this section we summarize our main results
on~\eqref{eq} in terms of $k$. We also give an account of
 the tools that we have developed for this study that we believe
might be interesting by themselves. We start by introducing the notations and
definitions used in the paper.

\subsection{Notations and definitions}

Given a periodic sequence $\{a_n\}$ of primitive period $k$ we will
say that its {\it rank} is $m$ if
$$
\mathrm{Card}\{a_1,a_2,\ldots,a_k\}=m\in\N.
$$
The values $a_1,\ldots,a_k$ will be usually called {\it parameters}.
In our context the recurrence~\eqref{eq} is called {\it persistent}
if for any sequence $\{x_n\}$ there exist two real positive
constants $c$ and $C,$ which depend on the initial conditions, such
that for all $n,\, 0<c<x_n<C<\infty.$

For each $k$, the {\it composition maps} are
\begin{equation}\label{MAPLYN} F_{a_k,\ldots
,a_2,a_1}:=F_{a_k}\circ\cdots\circ F_{a_{2}}\circ
F_{a_1}\end{equation}
 where each $F_{a_i}$ is defined by
$$F_{a_i}(x,y)=\left(y,\frac{a_i+y}{x}\right)$$ and $a_1,a_2,\ldots,a_k$ are the $k$
elements of the cycle. When there is no confusion, for the sake of
shortness, we also will use the notation
  $
F_{[k]}:=F_{a_k,\ldots ,a_2,a_1}.
  $
Note that these maps are birational maps and are always well-defined
in the open invariant set $\Qm=\{(x,y):x>0,y>0\}\subset\R^2.$
Moreover
\[
(x_1,x_2)\xrightarrow{F_{a_1}}(x_2,x_3)\xrightarrow{F_{a_2}}(x_3,x_4)
\xrightarrow{F_{a_3}}(x_4,x_5)\xrightarrow{F_{a_4}}(x_5,x_6)\xrightarrow{F_{a_5}}\cdots
\]
and in general,
\[
F_{[k]}(x_1,x_{2})=(x_{k+1},x_{k+2}).
\]

There are two  concepts coexisting in this context, {\it the
non-autonomous invariants} and {\it the first integrals}, that we
will use in this paper. Given a difference equation of the
form~\eqref{eq}, a non-autonomous invariant is a function
$V(x,y,n)$, such that
$$
 V\left(x_{n+1},x_{n+2},n+1\right)=V(x_{n},x_{n+1},n),
$$
for all initial conditions and all $n\in\mathbb{N}.$ On the other
hand, when the difference equation  has $k$-periodic coefficients a
first integral is a function $H$, which is a first integral for the
discrete dynamical system generated by $F_{[k]},$ that is
$H(F_{[k]}(x,y))=H(x,y),$ for all points in an open set. In terms of
the recurrence
\[H(x_{n+k},x_{n+k+1})=H(x_n,x_{n+1}),\]
for all initial conditions $(x_{n},x_{n+1})$. We will relate both
concepts in Section~\ref{se:3}.

Two analytic functions
$P,Q:\mathcal{U}\subset\mathbb{C}^2\to\mathbb{C}$ are said to be
coprime if the points of  the set  $\{(x,y)\in \mathcal{U}\,:\,
P(x,y)=Q(x,y)=0\}$ are isolated. A function $H=P/Q,$ with $P$ and
$Q$ coprime, will be called a {\it meromorphic function}. A {\it
meromorphic first integral} of an analytic map $F:\mathcal{U}\to
\mathbb{C}^2$ is a meromorphic function $H=P/Q$ such that
\[P(F(x,y))Q(x,y)=P(x,y)Q(F(x,y))\quad \mbox{for all}\quad (x,y)\in\mathcal{U}.\]
Observe that from this definition $H(F(x,y))=H(x,y)$ for all points
of $\mathcal U$ for which both terms of this last equality are
well-defined. When $P$ and $Q$ are polynomials then it is said that
$H$ is a {\it rational first integral}. Similarly we can talk about
{\it meromorphic or rational invariants}, and in this sense we will
talk about \textit{rational} or \textit{meromorphic integrability}.

Finally, we will say that a planar map $F$ has {\it structurally
stable numerical chaos (SSNC)} when studying numerically several of
its orbits, we observe that it presents all the features of a
non-integrable perturbed twist map, that is: many invariant curves
and, between them, couples of orbits of $p$-periodic points (for
several values of $p$), half of them of elliptic type and the other
half of hyperbolic saddle type. Moreover the separatrices of these
hyperbolic saddles intersect transversally,  for instance see
\cite[Chapter 6]{AP}.

\subsection{Main results}\label{smr}

This subsection collects the general outlines of all our results
about the recurrence~\eqref{eq}, in terms of  $k.$ Figure~1 shows
some typical behaviors of the orbits of $F_{[k]}$. In fact we
consider some maps $G_{[k]}$, which are conjugate to $F_{[k]}$,
because the pictures are much more clear. See Lemma~\ref{gg} for the
definition of $G_{[k]}.$

\vspace{0.2cm}

\noindent{\bf Cases $\mathbf{k\in\{1,2,3,6\}}$ and other concrete
integrable cases}. For $k\in\{1,2,3\}$  it is already known that the
recurrences~\eqref{eq} are persistent. Moreover, either each
sequence $\{x_n\}$ is periodic, with period a multiple of $k$ or it
densely fills at most $k$ disjoint intervals of $\mathbb{R}^+$, see
\cite{CGM09}.  A key point for the proof is the existence of a
rational first integral for $F_{[k]}$. When $k=6$ we can also prove
the existence of a similar first integral (see Corollary~\ref{fi})
and the persistence of recurrence~\eqref{eq}. Moreover we are
confident that the same
 characterization of the sequences    $\{x_n\}$      holds but we have only been able
 to prove the result  when $F_{[6]}$ has a unique fixed point in the first quadrant,
 see Lemma~\ref{lk=5} and Proposition~\ref{c6}.

It is also satisfied that for any $k\ne5$ there are values
$a_1,\ldots,a_k$, with primitive period $k$ and high rank,
satisfying the property that all the sequences $\{x_n\}$ given by
\eqref{eq} are either periodic, with period a multiple of $k$, or
they densely fill at most $k$ disjoint intervals, see
Theorem~\ref{main}.

\vspace{0.2cm}

\begin{center}
\includegraphics[scale=0.20]{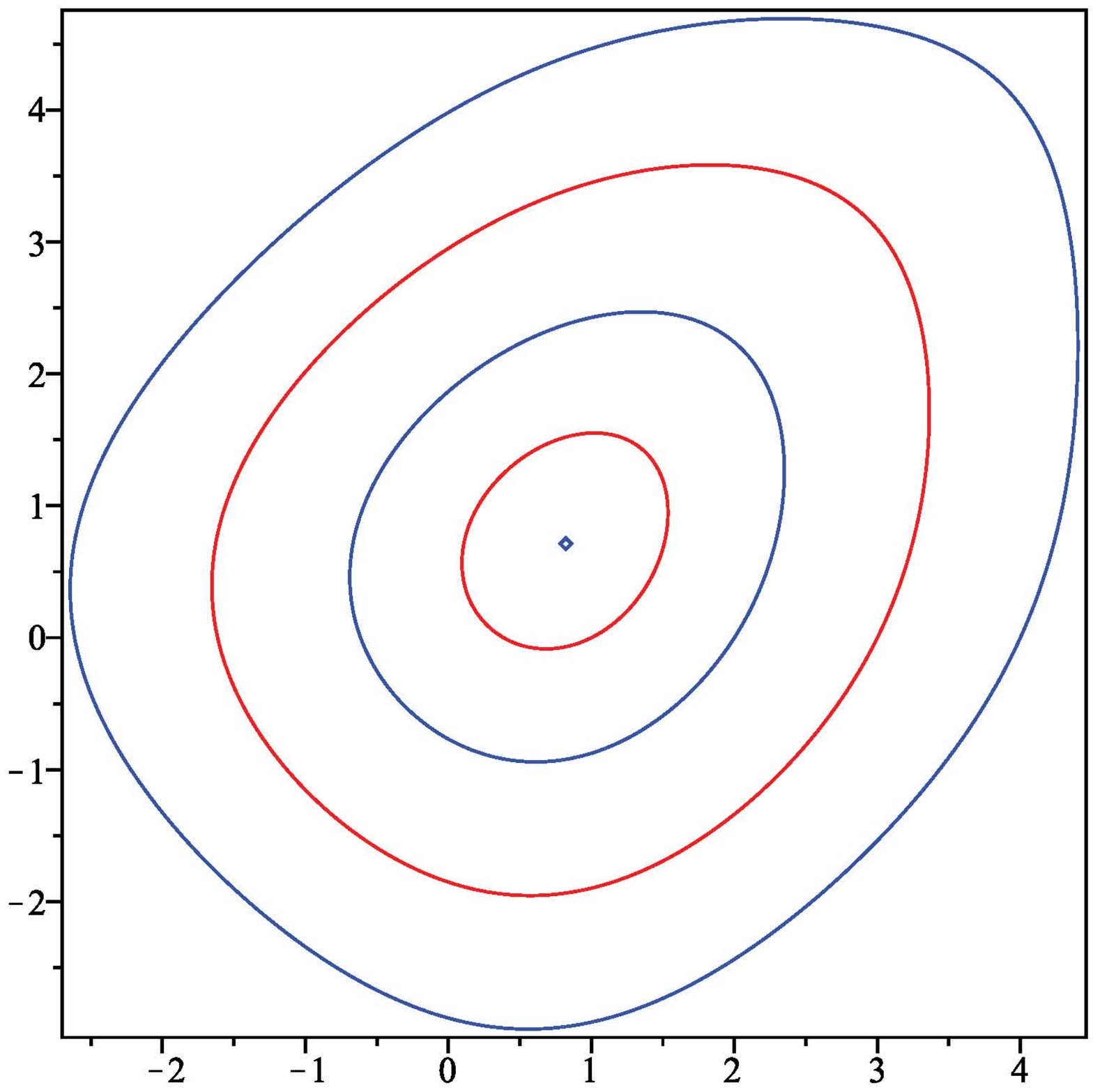}\hspace{1cm}
\includegraphics[scale=0.20]{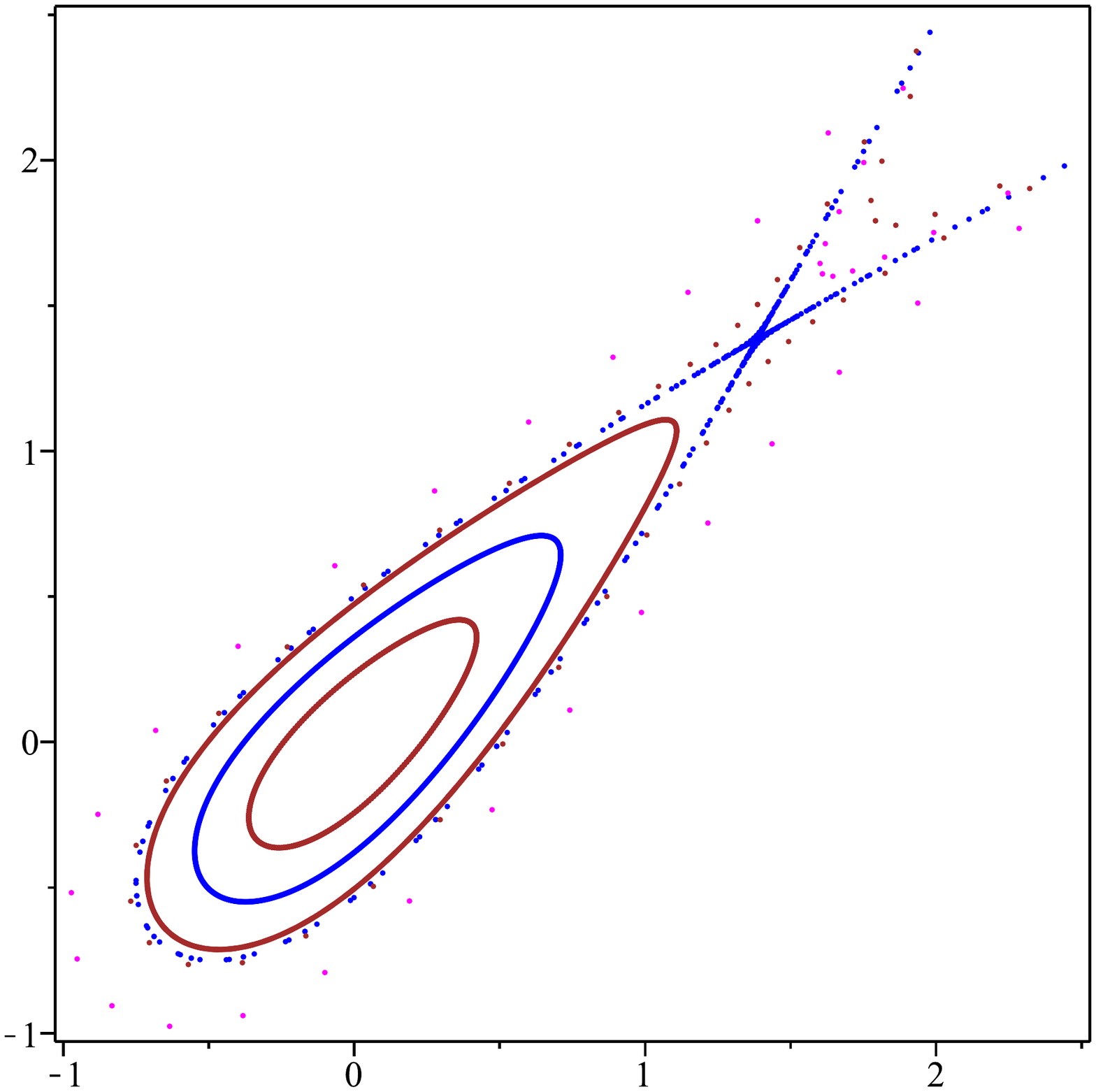}
\hspace{1cm}\includegraphics[scale=0.20]{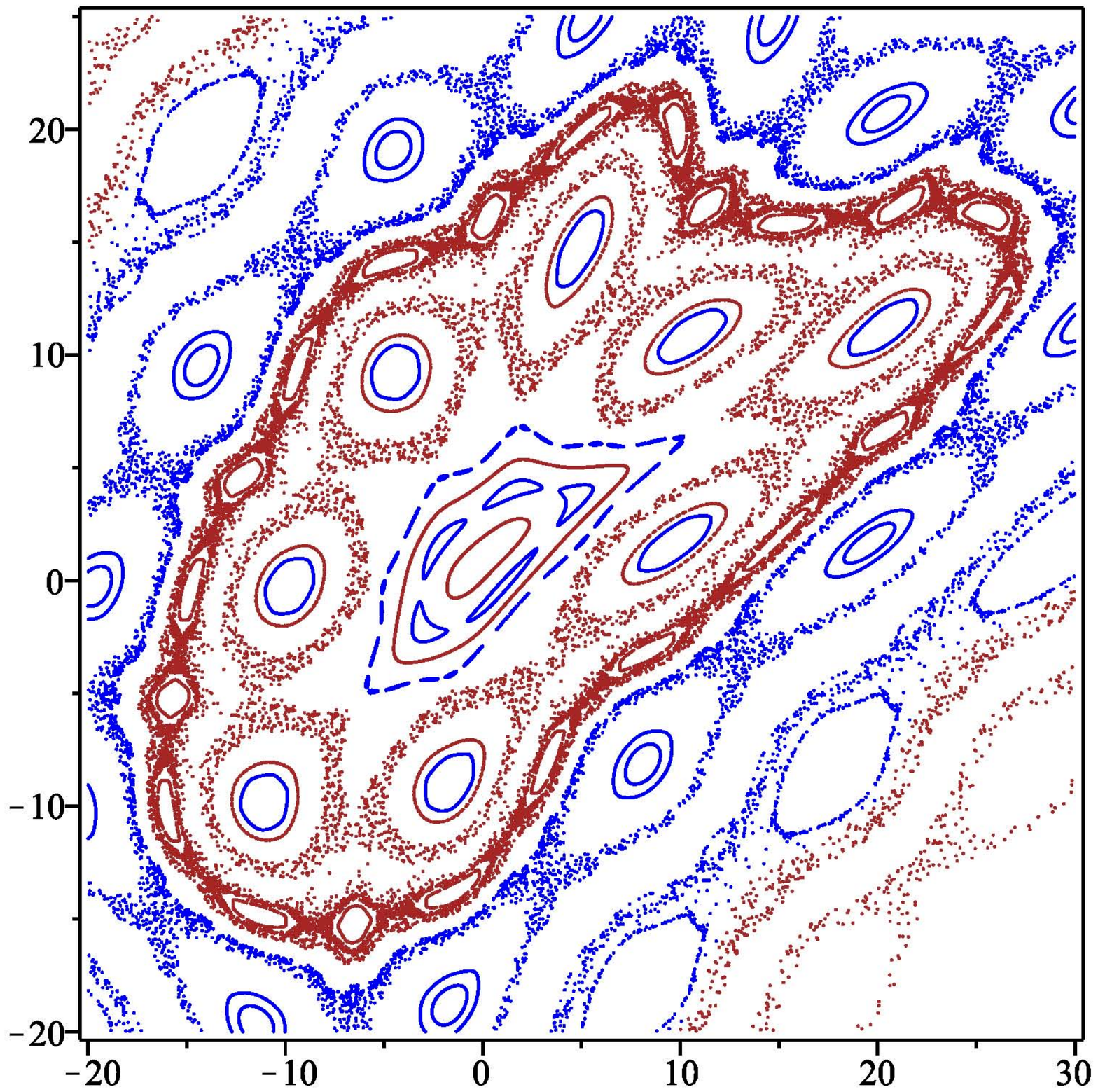}\\
\vspace{-1cm} \hspace{1cm} All $k$\hspace{4.5cm} $k=\dot{5}$\hspace{3.2cm}
$k\not\in\{1,2,3,5,6,10\}$

\end{center}
\begin{center}
Figure 1: Different possible behaviors of the orbits of $G_{[k]}$,
according to  $k$. Other behaviors are possible for $k$ being a
multiple of $5$.
\end{center}

\vspace{0.2cm}

\noindent{\bf Cases $\mathbf{k}$ being a multiple of $\mathbf{5}$}.
When $k$ is a multiple of 5 (from now on denoted by ${k=\dot{5}}$),
apart from the behaviors described above there appear others for an
open set of values of $a_j,j=1,\ldots,k$ and initial conditions.
 For instance we can
find sequences $\{x_n\}$ such that
\[
\liminf_{n\to\infty} x_n=0\quad\mbox{and}\quad \limsup_{n\to\infty}
x_n=+\infty,
\]
and others such that their adherence consists of $k$ points, see
Theorem~\ref{propo5}. The existence for $k=5$ of values of $a_j,
j=1,\ldots,5,$ for which  the sequence  $\{x_n\}$  has the first
 behavior has been already established in previous works, see
 \cite[Example 5.43.1]{CL} or \cite{deA},
  but only for very concrete initial conditions and  parameters $a_1,\ldots,a_5$.

Moreover in this case we can prove that for most values of the parameters the map
$F_{[k]}$ has no meromorphic first integral (see Theorem~\ref{nfi}). Furthermore
the phase portrait of the map $F_{[k]}$ does not always coincide with the ones
found in all the rational integrable cases, see for instance the second picture in
Figure 1. In this case, apart of the celebrated Lyness map $F_1$ which satisfies
$F_1^5=F_{1,1,1,1,1}=\operatorname{Id}$, there are values of the parameters
$a_1,\ldots,a_k$ such that $F_{[k]}^{m_k}=\operatorname{Id}$ (see
Corollary~\ref{gp}) and other for which the number of fixed points of the maps is
a 1-dimensional manifold, or $2,1,$ or $0$ points, see Lemma~\ref{lk=5}.  Finally,
when $k\ge15,$ there are  cases presenting  SSNC, see   Section~\ref{num}.
\vspace{0.2cm}

\noindent{\bf Cases $\mathbf{k\not\in\{1,2,3,5,6,10\}}$}. When $k\in\{4,7,11,15
\}$, for  some values of the parameters, $a_1,\ldots,a_k$, we have numerically
found SSNC, see  Section \ref{num}. In fact we prove in Lemma~\ref{chaos} that
based on these examples we can obtain values $a_1,\ldots,a_k$, all different, with
a similar behavior for all the remaining values of $k$. So for all these values of
$k$ there are situations for which the sequence $\{x_n\}$ can have different
behaviors to those given in the above situations. For instance there appear
sequences which fill more than $k$ intervals. Some concrete examples for $k=4$ are
shown in Section~\ref{num}.

Observe that as an application of our results we can show an interesting and
curious phenomena that can be understood as a kind of ``chaos regularization":
Consider a map $G=F_{[2]},$ which is a rationally integrable map, and a map
$H=F_{[4]}$ which has chaotic behavior, then both maps $G\circ H$ and $H\circ G$
are rationally integrable because are of type $F_{[6]}$. Thus $G$ regularizes $H$.

\smallskip

Finally note that the above results show that the only cases for
which recurrence \eqref{eq} can have a rational invariant for all
values of the parameters  $a_1,\ldots,a_k$ are $k\in\{1,2,3,6\}.$

\subsection{Main tools}\label{mt}

In this subsection we present   several  results that we have
obtained which we believe that are interesting by themselves. Other
technical results will be given in Section~\ref{pr}.

The first result is a necessary condition for the meromorphic
integrability of planar maps near a fixed point. Our approach
follows the guidelines of Poincar\'{e} when he  studied the same problem
for ordinary differential equations,  see \cite{LLZ} and the
references there in for the approach to ordinary differential
equations. In Section~\ref{k5}  we will apply the result below to
study the case $k=\dot{5}$.

\begin{teo}\label{th-int}
Let $F:\mathbb{C}^2\to\mathbb{C}^2$ an analytic map defined in
$\mathcal U$, an open neighborhood of the origin, such that
$F(0,0)=(0,0)$ and $DF(0,0)$ is diagonalizable with eigenvalues
$\lambda$ and $\mu$. Assume that $F$ has a meromorphic first
integral $H$ in $\mathcal{U}$.
 \begin{enumerate}[(i)]
 \item If  $\lambda\mu\ne0$ then there exists
 $(p,q)\in\mathbb{Z}^2,$ $(p,q)\ne(0,0)$,
 such that $\lambda^p\mu^{q}=1$.
  \item If $\lambda\ne0$ and $\mu=0$ then there exists  $n\in\mathbb{N}^+$
  such that $\lambda^n=1.$
  \end{enumerate}
\end{teo}

When the map $F:\mathcal{U}\subset\mathbb{R}^2\to \mathbb{R}^2$ is
real valued and of class $\mathcal{C}^2(\mathcal{U})$ the proof of
Theorem~\ref{th-int} can be adapted  following the same steps.
Taking into account that in this case, when $\lambda\in\mathbb{C}$
is an eigenvalue of $DF(0,0)$, then $\bar \lambda$ is also, and we
have to deal with the resonant condition
$\lambda^p{\bar\lambda}^q=1,$ we obtain the following result:

\begin{corol}Let $F:\mathcal{U}\subset\mathbb {R}^2\to\mathbb{R}^2$ be a
$\mathcal{C}^2(\mathcal{U})$ map
 such that $F(0,0)=(0,0)\in\mathcal{U}$ and $DF(0,0)$ is
diagonalizable, with eigenvalues $\lambda$ and $\mu$. Assume that
$F$ has a meromorphic  first integral $H$ in $\mathcal{U}$.
 \begin{enumerate}[(i)]
 \item If $\lambda,\mu\in\mathbb{R}$,  $\lambda\mu\ne0,$ then there exists $(0,0)\ne(p,q)\in\mathbb{Z}^2$
 such that $\lambda^p\mu^q=1$.
\item If $0\ne\lambda\in\mathbb{C}\setminus\mathbb{R}$ (hence $\mu=\bar \lambda$), then
either $|\lambda|=1$ or $\lambda=|\lambda|e^{i\theta}$ and there
exists $0\ne n\in\mathbb{N}$ such that $(e^{i\theta})^{2n}=1$.

  \item If $\lambda\ne0$ and $\mu=0$ then there exists a $n\in\mathbb{N}^+$
  such that $\lambda^n=1.$
  \end{enumerate}

\end{corol}

The above results will be applied to prove the meromorphic
non-integrability of many cases when $k=\dot{5},$  see
Theorem~\ref{nfi}. On the  other hand, next result will be the key
point to prove the existence of rational integrable cases for all
$k\ne5,$ see Theorem~\ref{main}.

For the recurrence (\ref{eq}) we look for non-autonomous invariants
of the form
\begin{equation}\label{formainvariant}
V(x,y,n)={\frac {\Phi_n(x,y)}{xy }}\end{equation} where
$$\begin{array}{rl}
         \Phi_n(x,y)=&A_{{n}}+B_{{n}}x+C_{{n}}y+D_{{n}}{x}^{2}+
         F_{{n}}{y}^{2}+G_{{n}}{x}^{3}+H_{{n}}{x}^{2}y+I_{{n}}x{y}^{2}\\
          &+
J_{{n}}{y}^{3}+K_{{n}}{x}^{4}
+L_{{n}}{x}^{3}y+M_{{n}}{x}^{2}{y}^{2}+N_
{{n}}x{y}^{3}+O_{{n}}{y}^{4},
        \end{array}
$$
with all the sequences of positive numbers. This method is
introduced in \cite{FJL} and the special form of $V$ is inspired by
this paper and the known invariant of the Lyness recurrences, see
\cite{BR,JKN,KN}. We prove:

\begin{teo}\label{noautonom}
If the  recurrence \eqref{eq} has an invariant of the form
\eqref{formainvariant} then $a_{n+6}=a_{n}$ and
\begin{align}\label{inv}
  \Phi_n(x,y)= & \,\,a_{n}F_{{n+1}}+(F_{{n+2}}+a_{{n+1}}F_{{n+1}})x+(F_{{n+1}}+
a_{{n}}F_{{n}})y+F_{{n-3}}{x}^{2}\nonumber\\
 & +F_{{n}}{y}^{2}+F_{{n-2}}{x}^{2}y +F_{{n-1}}x{y}^{2},
  \end{align}
where $\{F_n\}_n$ satisfies that $F_{n+6}=F_n$ and
$a_{n+1}F_{n+2}-a_nF_{n-3}=0.$
\end{teo}

\begin{corol}\label{fi} (i) The non-autonomous $k$-periodic
recurrence \eqref{eq} has invariants
of the form \eqref{formainvariant} if and only if $k\in\{1,2,3,6\}$.

(ii) The first integrals of the maps $F_{[k]}$, for
$k\in\{1,2,3,6\},$ corresponding to the invariants given in
Theorem~\ref{noautonom} are:
\begin{align*}
V_{a}(x,y)&=\frac{a+(a+1)x+(a+1)y+x^2+y^2+x^2y+xy^2}{xy},
\\
V_{b,a}(x,y)&= {\frac {ab +(a+b^2)x+(b+a^2)y
+bx^2+ay^2+ax^2y+bxy^2}{xy}},
\\
V_{c,b,a}(x,y)&= {\frac {ac+(a+bc)x+(c+ab)y+b{x}^{2}+b{y}^{2}+
c{x}^{2}y+ax{y}^{2}}{xy}},
\\
V_{f,e,d,c,b,a}(x,y)&={\frac {af+(a+bf)x+(f+ae)y+b{x}^{2}+e{y}^{2}+
c{x}^{2}y+dx{y}^{2}}{xy}}.
\end{align*}

\end{corol}

\begin{nota} Our  proof of Theorem~\ref{noautonom} does not require the sequence of parameters   $\{a_n\}$ to be  periodic.
\end{nota}

Observe that
\[
V_a(x,y)+2+a=\frac{(x+1)(y+1)(a+x+y)}{xy}
\]
is the usual first integral (invariant) of the map $F_a$ associated
to the classical  Lyness recurrence, see for instance \cite{BR}. It
is already known, see \cite{CGM09,JKN,KN}, that in the two and three
periodic cases the functions $V_{b,a}$ and $V_{c,b,a}$  are first
integrals of the maps $F_{b,a}$ and $F_{c,b,a},$ respectively. These
first integrals play a crucial role for the understanding of  the
recurrence \eqref{eq} when $k=2,3$, see again~\cite{CGM09}. To the
best of our knowledge the existence of a first integral  for the
general non-autonomous 6-periodic case was not known. In
Section~\ref{ri} we use it to describe the dynamics in this case.

 Also observe that, due to the form of the invariants, all the maps $F_{[k]}$,
for $k\in\{1,2,3,6\},$ preserve a foliation of the plane given by biquadratic
curves, which are elliptic except for a finite number of level sets. In fact these
maps are particular cases of the celebrated QRT family of planar maps. Perhaps a
further algebraic-geometric approach, like the one presented in \cite{D}, by
studying the maps induced by each $F_{[k]}$ on the corresponding elliptic surface
could give more information about the reason why the cases $k\in\{1,2,3,6\}$ are
special.

It is also interesting to notice that other integrable QRT maps with
periodic coefficients have been found recently \cite{RGW}, as well
as another major family of maps: the Hirota-Kimura-Yahagi type ones
\cite{GRT2011}.

As we have already commented, it is known that for very concrete
values of $a_1,\ldots,a_5$ and suitable initial conditions, the
behavior of $\{x_n\}$ is different to the ones appearing when
$k\in\{1,2,3\}$ and in particular \eqref{eq} is  non-persistent, see
 \cite[Example 5.43.1]{CL} or \cite{deA}. This behavior
can also be seen considering
$$F_{a,1,1,1,1}(x,y)=\left(x,{\frac { \left( x+a \right) y}{1+x}}\right).$$
Since  $F_{a,1,1,1,1}(1,y)=\left(1,\left( 1+a \right) y/2\right)$ it
is clear that for $a>1$ the orbits of the points of the form $(1,y)$
with $y\neq 0$ are unbounded.

Our next result allows to establish, when $k=\dot{5}$, the
non-persistence of the recurrence~\eqref{eq}  for many values
$a_1,\ldots,a_k.$

\begin{teo}\label{propo5} Consider recurrence~\eqref{eq} for $k=\dot{5}$.
Set
$$
\phi_{i}:=\prod\limits_{\tiny \begin{matrix} n\equiv
i\,(\operatorname{mod}\,5)\\ n=1,\ldots,k\end{matrix}} a_n, \,\mbox{
for }\quad i=1,2,\ldots, 5.
$$
If for all $i=1,\ldots,5$, $\phi_i\neq 1$ and
\begin{equation}\label{relaciophi}
\min_{i=1,\ldots, 5}\{\phi_i\}<1<\max_{i=1,\ldots, 5}\{\phi_i\},
\end{equation}
then~\eqref{eq} is non-persistent. In fact, for an open set of
initial conditions, $\liminf_{n\to\infty}(x_n)=0$ and
$\limsup_{n\to\infty}(x_n)=+\infty$.
\end{teo}

In the above result, the open set of initial conditions for which
the result holds is sometimes the whole first quadrant $\Qm$.
 For instance, this is the case when $k=5$ and $a_1=a, a_2=ac, a_3=c, a_4=1/a$ and
$a_5=1/(ac)$, when $a>1$ and $ac>1$, because
\begin{equation}\label{eqlin}
F_{\frac{1}{ac},\frac{1}{a},c,ac,a}(x,y)=\left(\frac{x}{a},\frac{y}{ac}\right),
\end{equation}
is a linear map with a stable node at the origin.

The rest of the paper is organized as follows. In Section~\ref{pr}
we introduce some preliminary results, while in Section~\ref{se:3}
we prove the main tools described in Section~\ref{mt}.
Section~\ref{ri} is devoted to the cases for which we find rational
integrability, while in Section~\ref{k5} we prove the
non-integrability results when $k$ is a multiple of five. Finally,
in Section~\ref{num}  we present some numerical evidence of chaos.

 \section{Preliminary results}\label{pr}

This section contains some technical preliminary results and other
known results that we will use in the proofs given in subsequent
sections.

The next result will be useful for our numerical simulations. As we
will see, some new variables allows to ``observe'' much better the
numerical non-integrability studied in Section~\ref{num}. Its proof
is straightforward.

\begin{lem}\label{gg} The sequence \eqref{eq} in the variables
$z_n:=\log (x_n)$ becomes
\[
z_{n+2}=-z_n+\log(a_n+\exp(z_{n+1})),
\]
and the corresponding maps $F_{[k]}$ are conjugate to $G_{[k]}$, where $
G_{[k]}=G_{a_k,a_{k-1}\ldots,a_2,a_1}, $
\[
G_{a}(x,y)=(y,-x+\log(a+\exp(y)),
\]
and each $G_a$ is  defined on the whole plane, $\R^2$.
\end{lem}

Observe that the maps  $G_{[k]}$ are area preserving. In fact, it is
easy to see that the maps $F_{[k]}$ and $G_{[k]}$ satisfy the
following properties:

\begin{lem}\label{mus} For every
choice of positive numbers  $a_1,\ldots,a_k,$
\begin{enumerate}[(i)]

\item The map $F_{[k]}$ preserves the measure $\textbf{m}(B)=\int
_B\frac1{xy}\,dx\,dy$ or, in other words, it preserves the
symplectic form $\omega:=\frac{1}{xy} dx\wedge dy$. In consequence,
it holds that $$\mu(F_{[k]}(x,y))=\det(DF_{[k]} (x,y))\,\mu(x,y),$$
where $\mu(x,y)=xy.$

\item The map $G_{[k]}$ preserves the Lebesgue measure $\textbf{n}(B)=\int
_B\,dx\,dy.$ That is, it preserves the canonical symplectic form
$dx\wedge dy$. Hence, it holds that $\det(DG_{[k]} (x,y))\equiv1$.
\end{enumerate}
\end{lem}

 A related issue in connection with the above
lemma  is the fact (\cite[Thm. 1]{B}) that  the group of symplectic birational
transformations of the plane (which is the group of birational transformations of
$\mathbb{C}^2$ which preserve the differential form $\omega$) is generated by
 compositions of the Lyness map $F_1$ (the
$5$-periodic case of $F_{[1]}$, with $a=1$); an scaling; and a map
of the form $(x, y)\to (x^ay^b, x^cy^d)$ where the matrix $\left(
                                                      \begin{array}{cc}
                                                        a & b \\
                                                        c & d \\
                                                      \end{array}
                                                    \right)
\in  SL(2,\mathbb{Z})$. This was conjectured by Usnich in \cite{U} and recently
proved by Blanc in the above mentioned reference.

The next result allows to know the dynamics of each $F_{[k]}$ when
the map has a smooth first integral, its applicability  to
characterize the dynamics of any integrable $F_{[k]}$ is guaranteed
by Lemma \ref{mus}.

\begin{teo}[\cite{CGM}]\label{teor}
Let $\U\subset \mathbb{R}^2$ be an open set and let $F:\U\rightarrow
\U$ be a diffeomorphism such that it has a smooth regular first
integral $V:\U\rightarrow \R$ and there exists a smooth function
$\mu:\U\rightarrow \R^+$ such that for any $(x,y)\in \U,$
$\mu(F(x,y))=\det(DF(x,y))\,\mu(x,y).$ Then the following holds:
\begin{enumerate}[(i)]
\item If a  level set $\Gamma_h:=\{(x,y)\in\U\,:\, V(x,y)=h\}$
is a simple closed curve invariant under~$F$, then the map $F$
restricted to  $\Gamma_h$ is conjugate to a rotation.

\item If $\Gamma_h$ is diffeomorphic to an open interval curve and invariant under
$F$, then the map~$F$ restricted to  $\Gamma_h$ is conjugate to a
translation.
\end{enumerate}
\end{teo}

The next result proves that all the fixed points $\mathbf{p}$ of
some map $F_{[k]},$ in $\R^2$ such that
$F_{a_j,a_{j-1},\ldots,a_1}(\mathbf{p})$, for all $j\le k$, is well
defined are resonant. As we will see in Proposition~\ref{final}, the
hypothesis on the maps $F_{a_j,a_{j-1},\ldots,a_1}(\mathbf{p})$ is
unavoidable because for $j=\dot{5}$ there appear some cancelations
that make that the maps $F_{[5j]}$ have as fixed point
$\mathbf{p}=(0,0)$ and this point can be of saddle type with
arbitrary eigenvalues. In fact this property will be the key point
for our proof of non-existence of meromorphic first integrals for
most $F_{[k]}, k=\dot{5}$, see Theorem~\ref{nfi}.

\begin{propo}\label{determinantesu} Let $\mathbf{p}\in\mathbb{R}^2$ be a fixed point of
a  composition map $F_{[k]}$ and such that
$F_{a_j,a_{j-1},\ldots,a_1}(\mathbf{p})$ is well defined  for all
$j\le k$.   Then
 \begin{equation}\label{det}
\det\,(DF_{[k]}(\mathbf{p}))=1.
\end{equation}
\end{propo}

\begin{proof}
From Lemma \ref{mus} we know that  $G_{[k]}$ is a symplectic map that preserves
the canonical form, hence $\det(DG_{[k]})=1 $. On the other hand, since $F_{[k]}$
and $G_{[k]}$ are conjugated maps (Lemma \ref{gg}), the Jacobian matrices of these
maps at their corresponding fixed points have the same eigenvalues and the result
follows.~\end{proof}

As a consequence of the above  Proposition we have:
\begin{corol}\label{corolofftopic}
Let $\mathbf{p}\in\mathbb{R}^2$ be a $m$-periodic point of a
composition map $F_{[k]}$  such that
$F_{a_j,a_{j-1},\ldots,a_1}(\mathbf{p})$ is well defined  for all
$j\le  k m$.
   Then
 \begin{equation*}\label{det2}
\det\,(DF_{[k]}^m(\mathbf{p}))=1.
\end{equation*}
\end{corol}

The following  lemma studies the number of fixed points of $F_{[k]}$
for $k=4,5,6.$

\begin{lem}\label{lk=5}
\begin{enumerate}[(i)]

\item There is a unique fixed point of $F_{d,c,b,a}$ in $\Qm$ and it
satisfies
\begin{equation*}
\begin{cases}
&x=y^2+(a-c)y-d,
\\
&y=x^2+(d-b)x-a.
\end{cases}
\end{equation*}

\item There  are either 0,1,2 or a continuum of fixed points of
$F_{e,d,c,b,a}$ in $\Qm$ and they satisfy
\begin{equation*}
\begin{cases}
&(b-1)x+(1-d)y+(a-e)=0,
 \vspace{0.1cm}\\
&cxy-(e+x)(a+y)+bx+y+a=0.
\end{cases}
\end{equation*}

\item Let $\mathcal{F}_{f,e,d,c,b,a}\subset\Qm$ be the set of fixed
points of map $F_{f,e,d,c,b,a}$ and let
$\mathcal{S}_{f,e,d,c,b,a}\subset\Qm$ be the set of singular points
of its first integral $V_{f,e,d,c,b,a}$ given in Corollary~\ref{fi}.
Then $\mathcal{F}_{f,e,d,c,b,a}=\mathcal{S}_{f,e,d,c,b,a}$ and both
sets coincide with the set of points of  $\Qm$ satisfying
\begin{equation*}
\begin{cases}
y^2=\dfrac{(f+x)(a+bx)}{e+dx},
 \vspace{0.1cm}\\
x^2=\dfrac{(a+y)(f+ey)}{b+cy}.
\end{cases}
\end{equation*}
Moreover $\operatorname{card}(\mathcal{F}_{f,e,d,c,b,a})\ge1.$

\end{enumerate}

\end{lem}

\begin{proof}
(i) The set of fixed points of  $F_{d,c,b,a}$ is exactly the set of
points satisfying
\[
F_d(F_c(F_b(F_a(x,y))))=(x,y),
\]
but it is not easy to handle these two equations. On the other hand,
the equivalent condition
\[
F_b(F_a(x,y))=F^{-1}_c(F^{-1}_d(x,y)),
\]
lead us to the system of the statement. Clearly both parabolas meet
at a unique point in~$\Qm$.

(ii) Studying the condition
\[
F_c(F_b(F_a(x,y)))=F^{-1}_d(F^{-1}_e(x,y))
\]
we obtain the system of the statement. The $x$-coordinate of a fixed
point has to satisfy the quadratic equation
\[
(c-1)(b-1)x^2+(2e-1+bd+ac-ec-eb-ad)x+(e-1)(e-ad)=0.
\]
From this equation we  easily obtain the result. Notice that simple
cases having infinitely many fixed points appear for instance when
$b=d=1$ and $e=a.$

(iii) The two conditions given by
\[
F_c(F_b(F_a(x,y)))=F^{-1}_d(F^{-1}_e(F^{-1}_f(x,y))),
\]
directly lead to the system of the statement. The set of singular
points of $V_{f,e,d,c,b,a}$ is formed by the points satisfying
\[
\{(x,y)\in\Qm\,:\, \frac{\partial}{\partial
x}V(x,y)=\frac{\partial}{\partial y}V(x,y)=0 \}.
\]
The two equations describing the above set exactly coincide again
with the two equations given in the statement. So
$\mathcal{F}_{f,e,d,c,b,a}=\mathcal{S}_{f,e,d,c,b,a}$. That
$\operatorname{card}(\mathcal{F}_{f,e,d,c,b,a})\ge1$ can be seen
studying the behavior of the functions, $\frac{(f+x)(a+bx)}{e+dx}$
and $\frac{(a+y)(f+ey)}{b+cy}$,
 near $0$ and $+\infty$.~\end{proof}

Finally, we will also use the following results:

\begin{lem}\label{lemanou}
The  map $F_{{1}/{a},c,ac,a}$ is
  conjugate to the Lyness' map $F_{{1}/{(ac^2)}}$.
\end{lem}

\begin{proof} Observe  that
$F_{\frac{1}{a},c,ac,a}(x,y)=\left(\dfrac{1+cx}{y},\dfrac{x}{a}\right)$
and
$F_{\frac{1}{ac^2}}(u,v)=\left(v,\dfrac{\frac{1}{ac^2}+v}{u}\right)$.
If we consider the  linear map
$$
\varphi(x,y)=\left(\frac{y}{c},\frac{x}{ac}\right),
$$
it holds that $F_{\frac{1}{ac^2}}=\varphi\circ
F_{\frac{1}{a},c,ac,a}\circ \varphi^{-1} $, as we wanted to
prove.\end{proof}

A nice consequence of the 5-global periodicity of the Lyness map
$F_1(x,y)=(y,(1+y)/x)$ and the above lemma is the following result:

\begin{corol}\label{gp} Recurrence~\eqref{eq} with $k=4$ and
$[a_1,a_2,\ldots]=[1/c^2,1/c,c,c^2,1/c^2,1/c,\ldots]$ is globally
20-periodic, i.e. $F^5_{c^2,c,1/c,1/c^2}(x,y)=(x,y)$ for all
$(x,y)\in\Qm.$
\end{corol}

\section{Proof of the main tools}\label{se:3}

This section is devoted to proving Theorems~\ref{th-int},
\ref{noautonom} and~\ref{propo5}.

\begin{proof}[Proof of Theorem~\ref{th-int}]
Write
\[
H(x,y)=\frac{P(x,y)}{Q(x,y)}=\frac{P_{\tilde n}(x,y)+O(\tilde
n+1)}{Q_{\tilde m} (x,y)+O(\tilde m+1)},
\]
where  $P_{\tilde n}$ and $Q_{\tilde m}$ are homogeneous polynomials
with degrees $\tilde n\ge0$ and $\tilde m\ge0$, respectively, and
$O(k)$ denotes terms of order at least $k$. Firstly we prove that it
is not restrictive to assume that $\tilde n\ge \tilde m$ and that if
$\tilde n=\tilde m$ then $P_{\tilde n}(x,y)/Q_{\tilde m}(x,y)$ is
not constant. Notice that if $H$ is a first integral then $1/H$ is
also. Hence we can assume that $\tilde n\ge \tilde m$. If $\tilde n=
\tilde m$ and $P_{\tilde n}=\eta Q_{\tilde n}$ for some
$0\ne\eta\in\mathbb{C}$, take $\tilde H=H-\eta$. Clearly $\tilde H$
is a new first integral of the form $\tilde H(x,y)=(O(\tilde
n+1))/(Q_{\tilde m}(x,y)+O(\tilde m+1)),$ as we wanted to see.

It is also clear that in a  neighborhood of the origin  we can
assume that $F(x,y)=(\lambda x+O(2),\mu y+O(2)).$

(i) We start studying the case $\lambda\mu\ne0$. By imposing that
$H$ is a first integral of $F$ in $\mathcal U$ we have that
\[
{P(F(x,y))}{Q(x,y)}={P(x,y)}{Q(F(x,y))}.
\]
 By taking the lower order
terms of the above equality we obtain
\begin{equation}\label{eqh}
{P_{\tilde n}(\lambda x,\mu y)}{Q_{\tilde m}(x,y)}= {P_{\tilde
n}(x,y)}{Q_{\tilde m}(\lambda x,\mu y)}.
\end{equation}
Define
\begin{equation*}
P_n(x,y)=\frac{ P_{\tilde n}(x,y)}{\gcd(P_{\tilde n}(x,y) ,Q_{\tilde
m}(x,y))},\quad Q_m(x,y)=\frac{ Q_{\tilde m}(x,y)}{\gcd(P_{\tilde
n}(x,y) ,Q_{\tilde m}(x,y))},
\end{equation*}
where $n\ge m$ are suitable non-negative integers. By using the
homogeneity of $P_n$ and $Q_m$,  equation \eqref{eqh} becomes
\begin{equation}\label{eqh2}
{\mu^n y^n P_n(\lambda x/(\mu y),1)}{y^m Q_m(x/y,1)} ={y^n
P_n(x/y,1)}{\mu^m y^m Q_m(\lambda x/(\mu y),1)},
\end{equation}
where notice that we have canceled the common factor of $P_{\tilde
n}$ and $Q_{\tilde m}$. By introducing the polynomials in one
variable $p_n(w)=P_n(w,1),$ $q_m(w)=Q_m(w,1),$ with respective
maximum degrees $n$ and $m$, and $\rho=\lambda/\mu,$ $w=x/y,$
equation \eqref{eqh2} becomes
\begin{equation}\label{eqh3}
{\mu^{n-m}p_n(\rho w)}{q_m(w)}={p_n(w)}{q_m(\rho w)},
\end{equation}
where we know that $p_n$ and $q_m$ are not identically zero and have
no common root.

Notice that equality \eqref{eqh3} implies that if $w=w^*$ is a root
of $p_n$ then $\rho w^*$ is also, and hence $\rho^l w^*$ for any
$l\in\mathbb{N}$, is a root of $p_n$. Since $p_n$ has at most $n$
roots, if $w^*\ne0$ we have that $\rho^k=1$ for some $k\le n,$
proving the theorem, because $\lambda^k\mu^{-k}=1$. A similar
reasoning can be done for $q_m$. Hence it only remains to study the
cases
\[
p_n(w)=aw^{\hat n},\quad 0\le \hat n \le n,\quad \mbox{and}\quad
q_m(w)=bw^{\hat m},\quad 0\le \hat m \le m,
\]
for some complex numbers $a$ and $b$, $ab\ne0.$ Remember that we
know that both polynomials  have no common roots. So, at least one
of the two numbers $\hat n$ or $\hat m$ has to be zero. In any case,
the equation \eqref{eqh3} becomes
\[
{\mu^{n-m}a\rho^{\hat n} w^{\hat n}{b} w^{\hat m}={a w^{\hat
n}}b\rho^{\hat m} w^{\hat m}},
\]
giving $\mu^{n-m}(\lambda/\mu)^{\hat n-\hat m}=\lambda^{\hat n-\hat
m}\mu^{n+\hat m-m-\hat n}=1,$ as we wanted to prove.

\vspace{0.2cm}

(ii) When $\lambda\ne0$ and  $\mu=0$  equation \eqref{eqh}, after
dropping the common factor of $P_{\tilde n}$ and $Q_{\tilde m}$,
becomes
\begin{equation*}
{P_n(\lambda x,0)}{Q_m(x,y)}={P_n(x,y)}{Q_m(\lambda x,0)}.
\end{equation*}
Notice that $P_n(x,0)=a x^n$, $Q_m(x,0)=bx^m$, and $(a,b)\ne(0,0),$
because, otherwise $P_n$ and $Q_m$ would have $y$ as a common
factor. Hence
\begin{equation*}
a\lambda ^n x^n{Q_m(x,y)}=b\lambda^mx^m{P_n(x,y)},
\end{equation*}
and so $ab\ne0.$ Therefore
$P_n(x,y)=a\lambda^{n-m}x^{n-m}Q_m(x,y)/b$. Since $P_n$ and $Q_m$
have no common factor, we get that $m=0$ and so $Q_m=b.$ Hence this
last equality becomes $P_n(x,y)=a\lambda^n x^n.$ Then $a
x^n=P(x,0)=a\lambda^n x^n,$ giving $\lambda^n=1,$ as we wanted to
prove. ~\end{proof}

To illustrate the above result in the next remark we present some
examples of maps having (or not having) meromorphic first integrals.

\begin{nota}
(i) The linear maps $F(x,y)=(\lambda x,\mu y)$, with $\lambda$ and
$\mu$ satisfying the resonant condition $\lambda^p\mu^{q}=1$ are the
simplest maps with meromorphic first integrals $H(x,y)=x^py^q.$

(ii) The map $F(x,y)=(x+y(x-y),0)$, with an eigenvalue 0,  has the
first integral $H(x,y)=(x-y+1)(y+1).$

(iii) For maps with identically zero linear part we can have
existence or not of meromorphic first integrals.  For instance the
map $F(x,y)=(x^2,xy)$ has the first integral $H(x,y)=x/y.$ On the
other hand, by using the same tools that in our proof of
Theorem~\ref{th-int}, we can prove that the
 map $F(x,y)=(x^2,y^2)$ has no meromorphic first integral.
\end{nota}

\begin{proof}[Proof of Theorem~\ref{noautonom}]
The condition that a function $V(x,y,n)$ of the form
(\ref{formainvariant}) is a non--autonomous invariant of the
recurrence (\ref{eq}) becomes
$$
 V\left(y,\frac{a_n+y}{x},n+1\right)-V(x,y,n)=0,
$$
for all $(x,y)\in\Qm$ and all $n\in\mathbb{N}$. Imposing that the
each one of the coefficients of the 31 monomials $x^iy^j$ vanishes
identically and playing a little bit with these conditions we obtain
that \eqref{inv} holds and moreover that
\begin{eqnarray}
a_{n+1}F_{n+2}-a_{n}F_{n-3}&=&0,\label{eee}\\
 F_{n+3}-F_{n-3}+a_{n+2}F_{n+2}-a_nF_{n-2}&=&0\nonumber.
\end{eqnarray}
From the first equation we get that
 \[
a_{n}=\frac{F_{n+2}}{F_{n-3}}\,a_{n+1}=\frac{F_{n+3}}{F_{n-2}}
\frac{F_{n+2}}{F_{n-3}}\,a_{n+2}.
 \]
Plugging this equation in the second one we obtain that
\[
\frac{F_{n+3}-F_{n-3}}{F_{n-3}}\left(
F_{n-3}-a_{n+2}F_{n+2}\right)=0.
\]
Clearly the above equation holds if either $\{F_n\}_n$ is a
6-periodic sequence or  $F_{n-3}=a_{n+2}F_{n+2}.$

In the first situation let us prove that if $\{F_n\}_n$ is a
$p$-periodic sequence, $p\in\{1,2,3,6\}$ then $a_n$ also has to be
$p$-periodic. Assume for instance that $p=3,$ then using~\eqref{eee}
we have
\begin{align*}
\frac{a_{n+3}}{a_{n}}=\frac{a_{n+3}}{a_{n+2}}\frac{a_{n+2}}{a_{n+1}}\frac{a_{n+1}}{a_{n}}=
\frac{F_{n-1}}{F_{n+4}}\frac{F_{n-2}}{F_{n+3}}\frac{F_{n-3}}{F_{n+2}}
=\frac{F_{n+2}}{F_{n+4}}\frac{F_{n+4}}{F_{n+3}}\frac{F_{n+3}}{F_{n+2}}=1,
\end{align*}
as we wanted to see. The other cases follow similarly.

In the second situation we have that $F_{n-3}=a_{n+2}F_{n+2}.$ Using
this equation and equality \eqref{eee} we have that
\[
a_{n+2}=\frac{F_{n-3}}{F_{n+2}}=\frac{a_{n+1}}{a_n}.
\]
which is a well-known 6-periodic recurrence, as we wanted to see.
Finally, using  \eqref{eee} six times we get that
$F_{n+6}=F_{n}$.~\end{proof}

Before proving Corollary~\ref{fi}, we explain here how
non-autonomous invariants and first integrals are related.  Consider
a recurrence with $k$-periodic coefficients and having a
non-autonomous invariant $V(x,y,n)$ that satisfies
$V(x,y,n)=V(x,y,n+k)$. Then $H(x,y):=V(x,y,1)$ is a first integral
of $F_{[k]}$. Conversely, if $H(x,y)$ is a first integral then
\[
V(x,y,n):=H(F_{a_k,a_{k-1},...,a_\ell}(x,y)),\quad\mbox{where}\quad
1\le\ell\le k,\quad n-\ell=\dot {k}
\]
is a non-autonomous periodic invariant of the recurrence. These relations are used
in the next corollary for constructing the first integrals of $F_{[k]},$
$k=1,2,3,6$ using the invariant found in the above theorem. Notice also that there
is another way to relate both concepts. Indeed, it is possible to replace the
non-autonomous $k$-periodic recurrence by an autonomous map on an enlarged phase
space of dimension $\mathbb{R}^2 \times \mathbb{R}^{k}=\mathbb{R}^{k+2}$, simply
considering the map $ (x_n, x_{n+1},a_1, a_2, ..., a_{k-1}, a_k) \to (x_{n+1},
x_{n+2},a_2, a_3, ..., a_k, a_1). $
 In this case, any  non-autonomous invariant is
just an ordinary first integral for the above map. This approach is not used in
this paper.

\begin{proof}[Proof of Corollary~\ref{fi}] (i) This result is proved
along the proof of Theorem~\ref{noautonom}, above.

(ii) We only give the details for  $k=6$. The other cases follow
similarly. We introduce the following notations for the
non-autonomous 6-periodic recurrences:
\begin{align*}
\{a_n\}=a_1,a_2,\ldots=&\,a,b,c,d,e,f,a,b,c,d,e,f,a,\ldots\\
\{F_n\}=F_1,F_2,\ldots=&\,1,\ell,m,n,o,p,1,\ell,m,n,o,p,1,\ldots
\end{align*}
From the relations $a_{n+1}F_{n+2}-a_nF_{n-3}=0,$ we obtain that
\[
\ell=\frac f e,\quad m=\frac a e, \quad n=\frac b e,\quad o=\frac c
e \quad \mbox{ and  }\quad p=\frac d e.
\]
Hence, with the notations of Theorem~\ref{noautonom} we get that
\begin{align*}
\Phi_1(x,y)&=a_1F_2+(F_3+a_2F_2)x+(F_2+a_1F_1)y+F_{-2}x^2+F_{1}y^2+F_{-1}x^2y+F_0
xy^2\\
&=a\ell+(m+b\ell)x+(\ell+a)y+nx^2+y^2+ox^2y+pxy^2\\
&=\dfrac{af+(a+bf)x+(f+ae)y+bx^2+ey^2+cx^2y+dxy^2}{e}.
\end{align*}
Hence $V_{f,e,d,c,b,a}(x,y)=e\Phi_1(x,y)/(xy)$ is a first integral
of $F_{[6]}$, as we wanted to prove.
\end{proof}

We first prove Theorem~\ref{propo5} for $k=5.$

\begin{propo}\label{lema5} Consider recurrence~\eqref{eq} with $k=5$  and $a_i\neq
1,$ for $i=1,2,\ldots,5$ and satisfying
\begin{equation}\label{minmax}
\min\{a_1,a_2,a_3,a_4,a_5\}<1<\max\{a_1,a_2,a_3,a_4,a_5\}.
\end{equation}
Then the recurrence~\eqref{eq} is non-persistent. Moreover for an
open set of initial conditions $\liminf_{n\to\infty}x_n=0$ and
$\limsup_{n\to\infty}x_n=+\infty$.
\end{propo}

\begin{proof}  A computation shows that $F_{[5]}(x,y)=\left(P_{1}(x,y),P_{2}(x,y)\right)$
where
$$
\begin{array}{l}
P_{1}(x,y)=\displaystyle{ \frac{x \left(
a_3xy+a_4{y}^{2}+a_2x+(a_1a_4+1)y+a_1 \right) }{ \left( a_1+y
 \right)  \left( a_1+a_2x+y \right)}},\vspace{0.1cm}\\
P_{2}(x,y)=\displaystyle{
 \frac{yN(x,y)}{ \left( a_1+a_2x+y+a_3xy \right)
\left(a_1+a_2x+y
 \right)}}
\end{array}
$$
and
\begin{align*}
N(x,y)=
&a_3{x}^{2}y+a_4x{y}^{2}+a_2{x}^{2}+(a_1a_4+a_2a_5+1)xy+a_5{y}^{2}+\\
&a_1(1+a_2a_5)x+2\, a_1a_5y+{a_1}^{2}a_5.
\end{align*}
First observe that, contrary to what it happens for $F_{[k]}, k<4$,
$F_{[5]}$ can be extended to a neighborhood of $\Qm$. Note also that
$(0,0)$ is a fixed point and $$ DF_{[5]}(0,0)=\left(
\begin{array}{cc}
                     \frac{1}{a_1} & 0 \\
                     0 & a_5
                   \end{array}
\right).
$$
So, under our hypotheses, the origin is a hyperbolic fixed point of
$F_{[5]}$. Arguing similarly with the shifted maps
$F_{a_1,a_5,a_4,a_3,a_2}$, $F_{a_2,a_1,a_5,a_4,a_3}$,
$F_{a_3,a_2,a_1,a_5,a_4}$, and $F_{a_4,a_3,a_2,a_1,a_5}$, we obtain
that the origin is also a hyperbolic fixed point for these maps,
with Jacobian matrices
$$
\left(
  \begin{array}{cc}
    \frac{1}{a_2} & 0 \\
    0 & a_1 \\
  \end{array}
\right),\, \left(
  \begin{array}{cc}
    \frac{1}{a_3} & 0 \\
    0 & a_2 \\
  \end{array}
\right),\, \left(
  \begin{array}{cc}
    \frac{1}{a_4} & 0 \\
    0 & a_3 \\
  \end{array}
\right)\,\mbox{ and }\, \left(
  \begin{array}{cc}
    \frac{1}{a_5} & 0 \\
    0 & a_4 \\
  \end{array}
\right),
$$
respectively.

Condition (\ref{minmax}) implies that at least there exist a
parameter with a value less than one, and other grater that one.
Since there are no parameters with value equal to one and the
sequence is cyclic we can choose two contiguous parameters and such
that $a_i<1$ and $a_{i+1}>1$. This implies that the origin is an
attractive fixed point for some of the five shifted maps. For
example, suppose that $a_2<1$ and $a_3>1$, then the origin is an
stable node for $F_{a_2,a_1,a_5,a_4,a_3}$.

Taking an initial condition $(x_0,y_0)$, with positive coordinates
and in the basin of attraction of the origin for the corresponding
shifted map we obtain that $\liminf_{n\to\infty}x_n =0$ for the
solution of equation (\ref{eq}) with initial condition $x_1=x_0$ and
$x_2=y_0$.

Recall   $F_{a_j}(x,y)=(y,(a_j+y)/x),$ $a_{j}\ne0$. Thus the fact
that $\limsup_{n\to\infty}x_n=+\infty$ follows because if some
$\{(x_{n_s},y_{n_s})\}_{n_s}$ tends to $(0,0)$ then the second
component of $\{F_{a_j}(x_{n_s},y_{n_s})\}_{n_s}$ tends to
$+\infty$.
\end{proof}

\begin{proof}[Proof of Theorem~\ref{propo5}]  Set $F_{[k]}$ for $k=5m$. We can write $
F_{[k]}=F_{a_{5m},\ldots,a_{5m-4}}\circ\ldots\circ
F_{a_{5},\ldots,a_{1}}, $ so  $F_{[k]}$ can be extended to a
neighborhood of $\Qm$. Furthermore, observe that
$$
DF_{[k]}(0,0)=DF_{a_{5m},\ldots,a_{5m-4}}(0,0)\circ\ldots\circ
DF_{a_{5},\ldots,a_{1}}(0,0)=\left(
\begin{array}{cc}
                     \frac{1}{\phi_2} & 0 \\
                     0 & \phi_1
                   \end{array}
\right).
$$
Similarly the Jacobian matrices of the shifted maps have the form
$$
\left(
  \begin{array}{cc}
    \frac{1}{\phi_{i+1}} & 0 \\
    0 & \phi_i \\
  \end{array}
\right).
$$

Arguing as in Proposition~\ref{lema5}, the relation
(\ref{relaciophi}) implies that at least there exists a couple
$(\phi_i,\phi_{i+1})$ such that one of the values is greater than
one and the other less than one. So the origin is an attractive
fixed point for some of the $m$ shifted maps.  Now the proof follows
again as  in Proposition~\ref{lema5}.\end{proof}

\section{Rational integrability and associated dynamics}\label{ri}

As we have already explained in Subsection~\ref{smr} the cases
$k=1,2,3$ are very similar and totally understood. From
Corollary~\ref{fi}  we can prove in the next proposition a similar
result when $k=6.$ Before stating the result we introduce the
following notation:
\begin{align*}
\mathbf{P}_{f,e,d,c,b,a}:=\{(a,b,c,d,e,f)\in(\R^+)^6\,:\,\mbox{
system } \eqref{sist} \mbox{ has a unique solution in }\Qm \}
\end{align*}

\begin{equation}\label{sist}
\begin{cases}
y^2=\dfrac{(f+x)(a+bx)}{dx+e},
 \vspace{0.1cm}\\
x^2=\dfrac{(a+y)(f+ey)}{b+cy}.
\end{cases}
\end{equation}

\begin{propo}\label{c6} For $k=6$
the recurrence  \eqref{eq} is persistent. Moreover if
$(f,e,d,c,b,a)\in\mathbf{P}_{f,e,d,c,b,a}$,  any  sequence $\{x_n\}$
generated by \eqref{eq} is either periodic, with period a multiple
of $6$, or it densely fills at most $6$ disjoint intervals of
$\mathbb{R}^+$.
\end{propo}

\begin{proof} We follow the same steps as in
the proof of \cite[Thm. 1]{CGM09}. To prove the persistence
of~\eqref{eq}, it suffices to show that each level curve
$\{(x,y)\,:\,V_{f,e,d,c,b,a}(x,y)=h\}\cap\Qm$ is bounded. Since
\[\frac{af}{xy}+\frac{a+bf}{y}+\frac{f+ae}x+\frac{bx}{y}+\frac{e{y}}{x}+
{c{x}}+dy=h,
\]
we know that
\[
\frac{f+ae}h \le x\le \frac h c\quad\mbox{and}\quad \frac{a+bf}h\le
y \le\frac h d
\]
and the persistence follows.

By Lemma~\ref{lk=5}.(iii), under our hypotheses, the set of fixed
points of $F_{[6]}$ and the set of singular points of $V_{[6]}$
coincide and consists of a single point. Following again the same
guidelines of the proof of \cite[Thm. 1]{CGM09}, which in turn is
based on \cite[Prop. 2.1]{BR2}, we prove that all the level curves
of $V_{[6]}$ in $\Qm$, apart from the fixed point, are diffeomorphic
to circles. Hence by using Lemma~\ref{mus} and
Theorem~\ref{teor}.(i) the proposition follows.
\end{proof}

\begin{nota}\label{rc6} We  believe  that
$\mathbf{P}_{f,e,d,c,b,a}$ is the whole of $(\R^+)^6$ but we have
not been able to prove this equality. In any case it is easy to find
sufficient conditions  to ensure that some $(a,b,c,d,e,f)$ belongs
to $\mathbf{P}_{f,e,d,c,b,a}.$ For instance, since
\begin{align*}
&\frac{\partial}{\partial x}\left(\dfrac{(f+x)(a+bx)}{e+dx} \right)=
\frac{bdx^2+2bex+ ae+ bef-adf}{(e+dx)^2},\vspace{2cm}\\
&\frac{\partial}{\partial y}\left(\dfrac{(a+y)(f+ey)}{b+cy} \right)=
\frac{cey^2+2bey+bf+abe-acf}{(b+cy)^2},
\end{align*}
when both numerators have no positive real roots the point
$(a,b,c,d,e,f)$ is in the set, because the functions that we have
derived are both increasing and so the curves defined by
system~\eqref{sist} cut in a single point.
\end{nota}

Next result collects our integrability results for any $k\ne5.$

\begin{teo}\label{main}
\begin{itemize}
  \item[(i)] For any $k\geq 15$, there exist sequences $\{a_n\}$ of prime
  period $k$ and rank $k$ such that
  $
F_{[k]}=F_{a_k,\ldots ,a_2,a_1}
  $ is rationally integrable and the corresponding recurrence~(\ref{eq}) is persistent.
  \item[(ii)] For any $k<15, k\ne 5,$ there exist sequences  $\{a_n\}$ of prime
  period $k$  with the ranks as in Table 1, such that
  $F_{[k]}$ is rationally integrable and the corresponding  recurrence~(\ref{eq})
  is persistent.

  \item[(iii)] Moreover it is possible to take
  in all the above cases  parameters $a_1,a_2,\ldots a_k$ such that
   each sequence $\{x_n\}$ is either periodic,
   with period a multiple of $k$, or it densely fills at most $k$
disjoint intervals of $\mathbb{R}^+$.
\begin{center}
\begin{tabular}{|c|c|c|c|c|c|c|c|c|c|c|c|}
  \hline
  $k$ & $1$ & $2$ & $3$ & $4$ & $5$ & $6$ & $7$ & $8$ & $9$ &$10\le k\le 14$&k$\ge15$ \\
  \hline
  $\operatorname{Rank}$ & $1$ & $2$ & $3$ & $4$ & {-} & $6$ & $3$ & $4$ & $5$ &$k-5$&$k$ \\
  \hline
\end{tabular}
\end{center}
 \centerline{Table 1. Possible ranks for integrable $F_{[k]}$.}

\end{itemize}
\end{teo}

\begin{proof}
 We start by
introducing some notation. Given $a_i$ and $c_i$ positive, we
consider the sets
$S_i:=\left\{\frac{1}{a_ic_i},\frac{1}{a_i},c_i,a_ic_i,a_i\right\}$.
Assume that $a_i$ and $c_i$ are such that
$\operatorname{Card}(S_i)=5$ and consider
$$
\Phi_i(x,y)=F_{\frac{1}{a_ic_i},\frac{1}{a_i},c_i,a_ic_i,a_i}(x,y)=
\left(\frac{x}{a_i},\frac{y}{a_ic_i}\right),
$$
where we have used expression~\eqref{eqlin}. Notice that $1\notin
S_i$. Given any natural number $m\ge1,$ we also consider $m$ sets
$S_1,S_2,\ldots,S_m$,  and define $\Phi^{[m]}=\Phi_m\circ
\Phi_{m-1}\cdots\circ \Phi_1$. Then
$$
\Phi^{[m]}(x,y)=\left(\frac{x}{\prod\limits_{i=1}^{m}a_i},\frac{y}{\prod\limits_{i=1}^m
a_i c_i}\right).
$$
When $m=0$ we consider $\Phi^{[{0}]}(x,y)=(x,y).$  Finally, choosing
the values of the parameters such that $\prod\limits_{i=1}^{m}a_i=1$
and $\prod\limits_{i=1}^{m}c_i=1$, we obtain that for all $m\ge3$,
$\Phi^{[m]}(x,y)=(x,y)$. Moreover, for these values of $m$, the
parameters $a_i$ and $c_i$ can be chosen such that
$\operatorname{Card}(\cup_{i=0}^m S_i)=5m$. Observe  also that  when
$m=1$ it is not possible to choose $a_1=1$ and $c_1=1$. When $m=2$
it is again possible but with $\operatorname{Card}(\cup_{i=0}^2
S_i)=5$.

Now we can start the proof of the theorem. First we study the cases
$k\le4$ and $k\ge 15.$ Consider $k=5m+\ell$ with
$\ell\in\{0,1,2,3,4\}$ and $m=0$ or $m\ge3$.

Now, taking
\begin{equation}\label{defih}\Psi(x,y)\,:=\,\left\{\begin{array}{lll}
 (x,y)&{\mbox{for}}&\ell=0,\\
F_a(x,y)&{\mbox{for}}&\ell=1,\\
F_{b,a}(x,y)&{\mbox{for}}&\ell=2,\\
F_{c,b,a}(x,y)&{\mbox{for}}&\ell=3,\\
F_{\frac{1}{a},c,ac,a}(x,y)&{\mbox{for}}&\ell=4,
\end{array}\right.\end{equation}
with suitable values of $a$, $b$ and $c$, we obtain that the orbits
of
$$
F_{[k]}(x,y):=\Psi\circ \Phi^{[m]}(x,y)=\Psi(x,y)
$$
are like the ones of $\Psi$ and the
$\operatorname{rank}(\{a_n\})=k.$ Then by using the known results
for $k=1,2,3$ and Lemma~\ref{lemanou} the result follows for
$k\ge15$ and $k\le4$.

When $k=6$ the result is proved in Proposition \ref{c6}. Finally for
$7\leq k\leq 14$ we consider
$$
\begin{array}{cll}
  k=7, & F_{b,a,1,1,1,1,1} & \mbox{ with Rank } 3, \\
  k=8, & F_{c,b,a,1,1,1,1,1} & \mbox{ with Rank } 4,  \\
  k=9, & F_{{1}/{c},c,ac,a,1,1,1,1,1} & \mbox{ with Rank } 5,
\end{array}
$$
and $\Psi\circ F_{\bar a\bar c,\bar a, 1/\bar c,1/(\bar a\bar
c),1/\bar a,1/(\bar a\bar c),1/\bar a,\bar c,\bar a\bar c,\bar
a}=\Psi$ for $10\leq k\leq 14,$ with $\bar a$ and $\bar c$ suitably
chosen. All these $F_{[k]}$ have rank $k-5$, as we wanted to prove.
\end{proof}

\section{Meromorphic non-integrability for the case $k=\dot{5}$}\label{k5}

Our main result is the following theorem:

\begin{teo}\label{nfi} For $k=\dot{5}$ and most values of $\{a_n\}$ the map
$F_{[k]}$  has no meromorphic first integral.
\end{teo}

Its proof is a consequence of the following result:

\begin{propo}\label{final} For $k=\dot{5}$ let $\phi_i,$
$i=1,\ldots,5$ be as in Theorem~\ref{propo5},
$$
\phi_{i}=\prod\limits_{\tiny \begin{matrix} n\equiv
i\,(\operatorname{mod}\,5)\\ n=1,\ldots,k\end{matrix}} a_n.
$$
Then if $
\{\phi_2,\phi_3,\phi_4,\phi_5\}\not\subset\{\phi_1^{r},\,r\in\Q\} $
the map  $F_{[k]}$ has no meromorphic first integral.
\end{propo}

\begin{proof} For simplicity, we prove result for the case $k=5,$
being the proof in the general case similar. Note that the above
condition reads as
\begin{equation}\label{nointe}
\{b,c,d,e\}\not\subset\{a^r,\, r\in\Q\}.
\end{equation}
If $F_{[5]}=F_{e,d,c,b,a}$ has a meromorphic first integral, the
same holds for all the other maps $F_{a,e,d,c,b}$, $F_{b,a,e,d,c}$,
$F_{c,b,a,e,d}$, and $F_{d,c,b,a,e}$. These five maps have the
origin $(0,0)$ as a fixed point and are  analytic in its
neighborhood (see the proof of Proposition~\ref{lema5}). Moreover
the corresponding couples of eigenvalues of their linear parts at
zero are $1/a,e;$ $1/b,a$; $1/c,b$; $1/d, c$ and $1/e,d$,
respectively. Hence, applying Theorem~\ref{th-int}, we obtain  the
following necessary conditions for the existence of a meromorphic
first integral
\[
a\,^{n_1}e\,^{m_1}=a\,^{n_2}b\,^{m_2}=b\,^{n_3}c\,^{m_3}=c\,^{n_4}d\,^{m_4}
=d\,^{n_5}e\,^{m_5}=1,\] for some $n_i,m_i\in\mathbb{Z}$,
$i=1,\ldots,5$. From these equalities we get that
$\{b,c,d,e\}\subset\{a^r,\, r\in\Q\}$. So the result follows.
\end{proof}

A simple corollary of the above result is:

\begin{corol} The map $F_{e,d,c,b,1}$ has a meromorphic first integral if and only if
$b=c=d=e=1.$
\end{corol}

Of course, all the known rationally integrable cases, like
$F_{[5]}=F_a^5,$ $F_{[10]}=F_{ba}^5,$ $F_{[15]}=F_{c,b,a}^3$ and
$F_{[20]}=F_{{1}/{a},c,ac,a}^5$ satisfy
 $\{\phi_2,\phi_3,\phi_4,\phi_5\}\subset\{\phi_1^{r},\,r\in\Q\}$.

\section{Numerical evidences of chaos}\label{num}

Our simulations show that there are examples of maps $F_{[k]}$ exhibiting SSNC
when $k\in\{4,7,8,11,15\}.$ These behaviors can be seen by plotting some orbits of
the corresponding conjugated maps $G_{[k]}$ for the $k$-periodic sequences of
parameters $2,2,\ldots,2,3$. More visual examples can be obtained by studying the
maps $G_{\delta,7,4,2}$, $G_{1,\delta,1,\delta,7,4,2}$,
$G_{\delta/2,1,\delta,1,\delta,7,4,2},$ and
$G_{1,\delta,2,\delta/2,1,\delta,1,\delta,7,4,2}$ and
$G_{7,4,2,\delta,1,\delta,2,\delta/2,1,\delta,\delta,\delta,7,4,2}$ with
$\delta=0.001.$ See some pictures in Figures~1 and~2.

\begin{figure}[h]
\begin{center}
\includegraphics[scale=0.30]{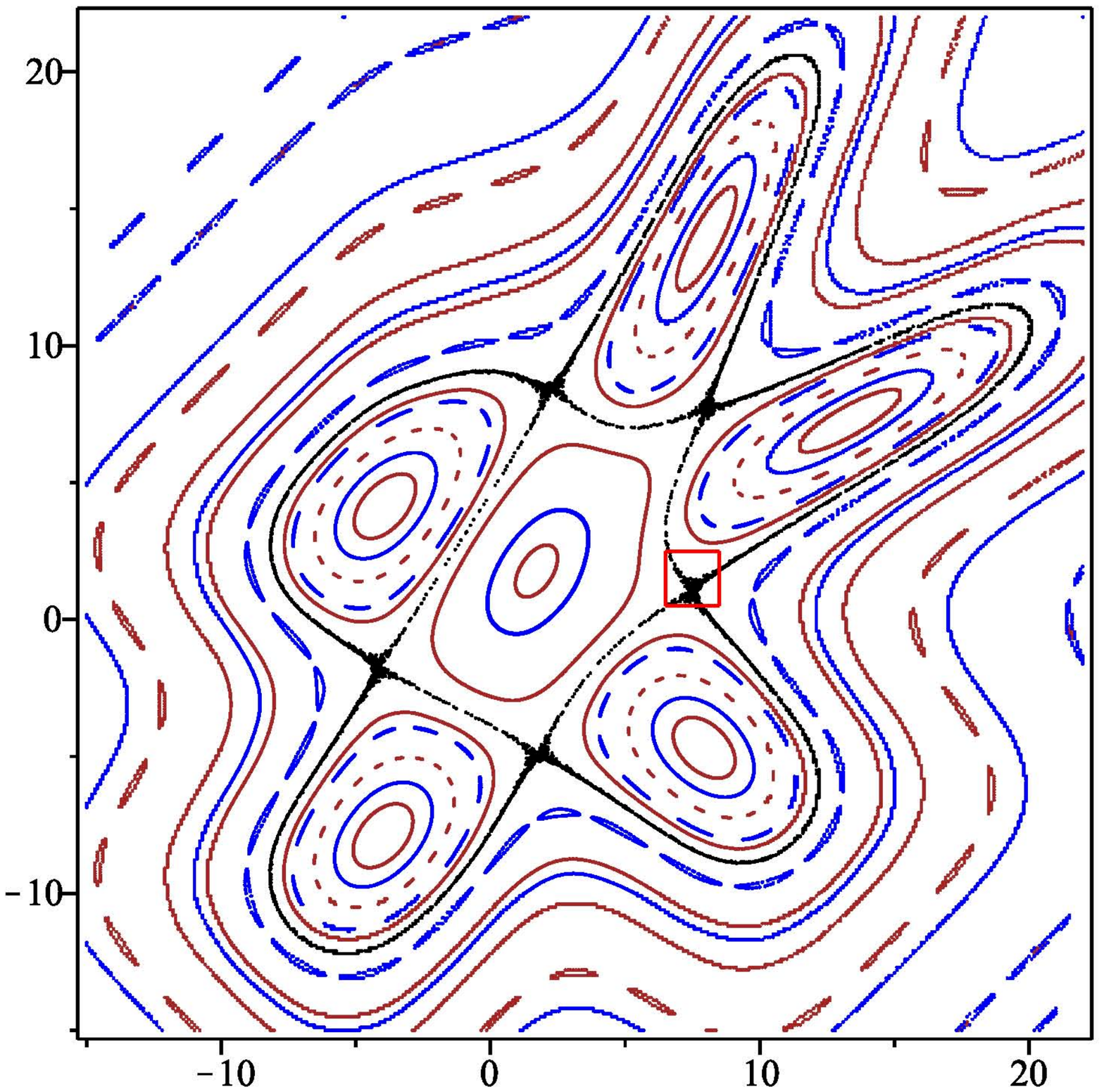}
\hspace{2cm}\includegraphics[scale=0.30]{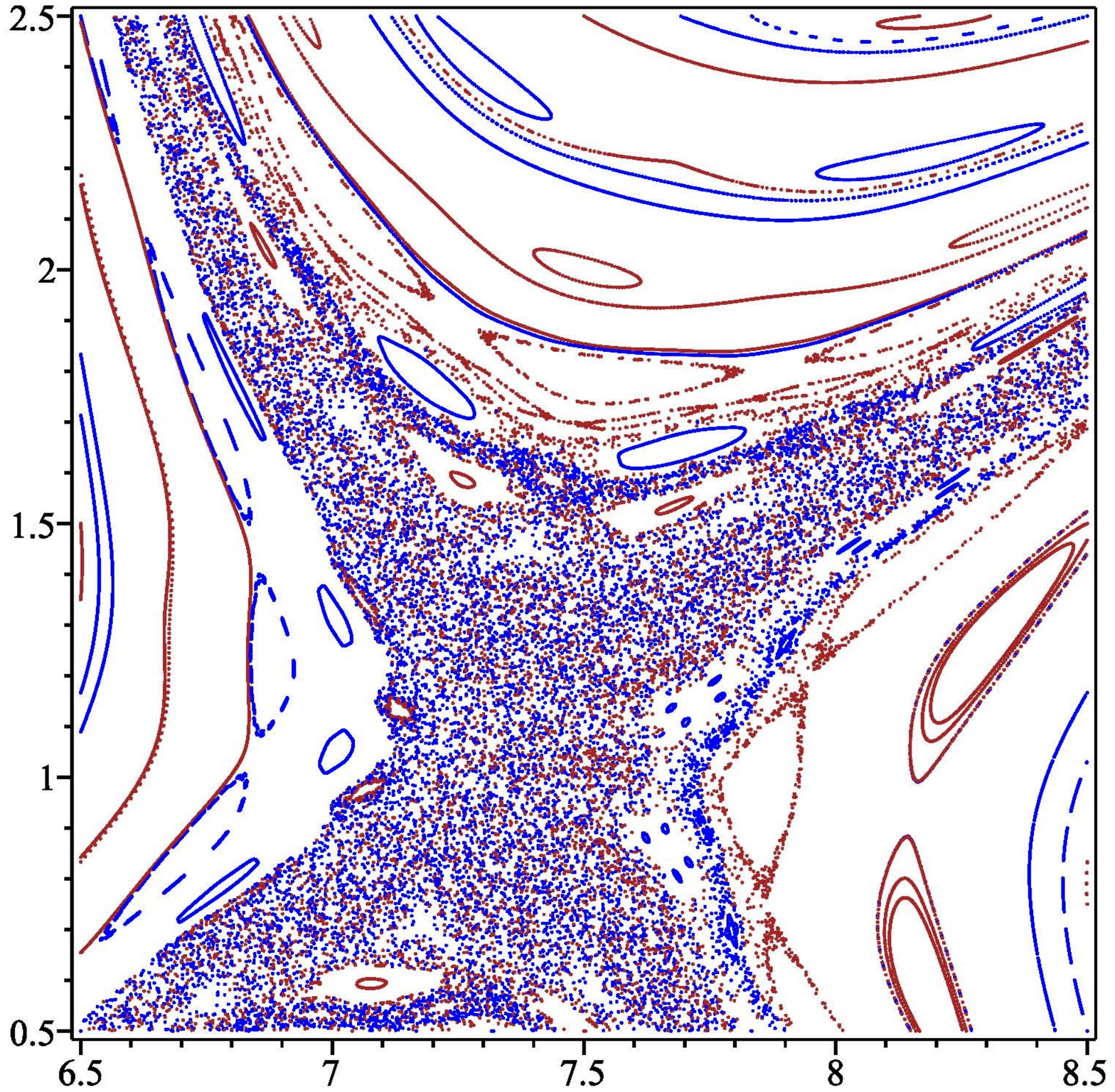}
\end{center}
\begin{center}
\vspace{-1.5cm} Figure 2: Some orbits of  $G_{0.001,7,4,2}$ and a
zoom with much more orbits.
\end{center}
\end{figure}

In fact we prove:

\begin{lem}\label{chaos} If there exist maps of the form $F_{[4]}$,$F_{[7]}$ and
 $F_{[11]}$ exhibiting SSNC then, for any $k\ge7,$ $k\not\in\{10,15\}$
 there exist maps of the form $F_{[k]}$ with
 $\operatorname{rank}(\{a_n\})=k$ having also SSNC.
\end{lem}

\begin{proof} For $m\ge0$, the
maps
\begin{align*}
F_{[5m+11]}&:=F^m_{1,1,1,1,1}\circ F_{[11]}=F_{[11]},\\
F_{[5m+7]}&:=F^m_{1,1,1,1,1}\circ F_{[7]}=F_{[7]},\\
F_{[5m+8]}&:=F^m_{1,1,1,1,1}\circ F_{[4]}\circ F_{[4]}=F_{[4]}\circ F_{[4]},\\
F_{[5m+4]}&:=F^m_{1,1,1,1,1}\circ F_{[4]}=F_{[4]},\\
F_{[5m+20]}&:=F^m_{1,1,1,1,1}\circ F_{[4]}^5= F_{[4]}^5,
\end{align*}
will also have SSNC. Note that, since $20\equiv0$, $11\equiv1$, $7\equiv 2$,
$8\equiv3$ and $4\equiv 4$ $(\operatorname{mod}\, 5)$,  the above maps cover all
the values of $k$ given in the statement. Since one
 of the features of these maps is the existence of transversal
 homoclinic points, which is a structurally stable property, we can
 perturb each of the corresponding $a_j$ by $a_j+\varepsilon_j$,
 with all the $\varepsilon_j$ sufficiently small, to obtain
 $k$-periodic sequences of parameters having SSNC and $\operatorname{rank}(\{a_n\})=k$, as
 we wanted to prove.
\end{proof}

We have (only numerically)  shown the existence of SSNC for $k\in\{4,7,11,15\}$,
but note that the above lemma allows to reduce all the other cases to these four
ones.

Although the complicated behavior of the maps, as in Figure 2, leads one to
believe that even there may be no upper bound for the number of intervals given by
the adherence of a sequence  (formed by the sequence itself and its accumulation
points), we only present here some simple examples. Concretely, for $k=4$ we give
a map $F_{[k]}$ and two sets of initial conditions such that the adherence of the
sequences $\{x_n\}$ generated by~\eqref{eq} consists of  more than $k$ intervals.

We have that  for $a=2,b=4,c=7$ and $d=0.001$,
\begin{itemize}
\item the sequence starting at $13.35, 7.27$ is formed by 20
intervals;

\item the sequence starting at $14.8,8.25$ is formed by 7 intervals.
\end{itemize}
For instance, this last assertion can be seen by making the phase
portraits of the orbit of $G_{0.01,7,4,2}$ starting at
$(14.8,8.25)$, which is formed by 5 islands, together with their
images through $G_{2}$, $G_{4,2}$ and $G_{7,4,2}$ and their
projections in the $x$-axis, see Figure~3. The property of the
existence of sequences $\{x_n\}$ generated by~\eqref{eq} such that
their adherence consists of  more that $k$ intervals should be true
for all the values of $k$ given in Lemma~\ref{chaos}. We also want
to comment that, for these values of $k$, the initial conditions
lying on the stable manifolds of the $q$-periodic saddle points of
$G_{[k]}$, also have a curious behavior: the adherence of $\{x_n\}$
is the sequence itself together with $q$ more points, corresponding
to the saddle points. In general $q$ also is greater than $k.$ On
the other hand, the most complicated orbits, that is the ones
between two big invariant curves, whose adherence seems to fill a
region of positive measure give rise to a single interval when we
consider their projections given by the sequence $\{x_n\}.$

\smallskip

\subsection*{A final remark} After the first version of this paper
was finished, one of the authors (A. Cima) and S. Zafar have
obtained a proof of the non rational integrability of generic maps
$F_{[4]}$  by computing its dynamical degree, \cite{CZ}.

\section*{Acknowledgements} We thank Joan Carles Tatjer for his
suggestion to introduce the maps $G_{[k]}$ for a better approach to
the numerical study of the maps $F_{[k]}$;  John A.G. Roberts who
gave us very useful information about the interest of classical
Lyness recurrences in the context of discrete integrability.


\begin{center}
   \includegraphics[scale=0.45]{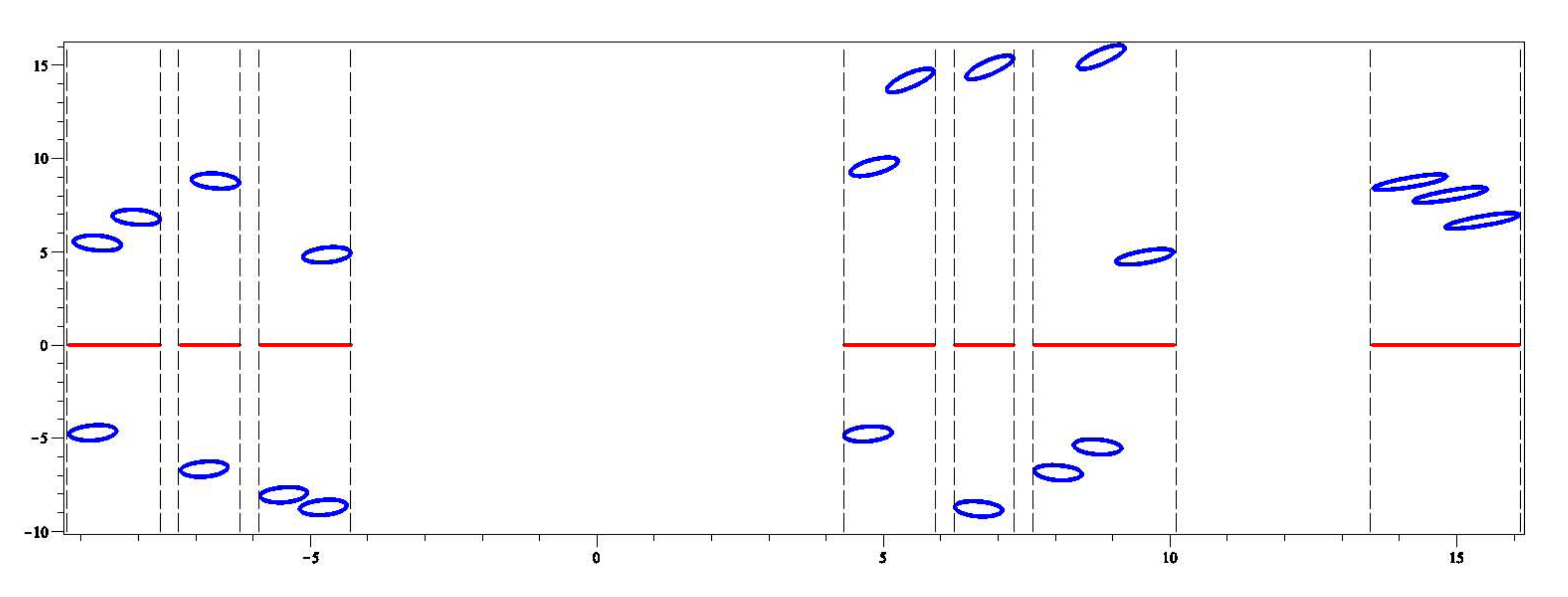}
\end{center}
   \begin{center}
   \noindent{ Figure 3. Case $k=4$. An orbit of $G_{[4]}$, their images through
   $G_{a_i,a_{i-1},\ldots,a_1},$ $i=1,2,3,$\\ and the projection
   corresponding to $\{x_n\}$.}
   \end{center}

\end{document}